\documentclass{amsart}

\usepackage{amsmath} 
\usepackage{amssymb}

\newtheorem{theorem}{Theorem}[section]

\newtheorem{corollary}[theorem]{Corollary} 

\theoremstyle{definition}
\newtheorem{definition}[theorem]{Definition}

\theoremstyle{remark}

\numberwithin{equation}{section}
\newcommand{\forces}{\Vdash} 
\newcommand{\bV}{{\bf V}} 

\newcount\skewfactor
\def\mathunderaccent#1#2 {\let\theaccent#1\skewfactor#2
\mathpalette\putaccentunder}
\def\putaccentunder#1#2{\oalign{$#1#2$\crcr\hidewidth
\vbox to.2ex{\hbox{$#1\skew\skewfactor\theaccent{}$}\vss}\hidewidth}}
\def\name{\mathunderaccent\tilde-3 }

\newcommand{\conc}{{}^\frown\!}
\newcommand{\lh}{{\rm lh}\/}
\newcommand{\rest}{{\restriction}}
\newcommand{\Dom}{{\rm Dom}} 
\newcommand{\Rng}{{\rm Rng}}
 
\newcommand{\suc}{{\rm succ}}
\newcommand{\vtl}{\vartriangleleft}
\newcommand{\baire}{{}^{\omega}\omega}
 
\newcommand{\bbD}{{\mathbb D}}

\newcommand{\cF}{{\mathcal F}}

\newcommand{\bbP}{{\mathbb P}}

\newcommand{\bbQ}{{\mathbb Q}}
\newcommand{\mbR}{{\mathbb R}}

\newcommand{\bbS}{{\mathbb S}}

\newcommand{\st}{{\bf st}} 
\newcommand{\FP}{{\rm FP}}
\newcommand{\tGs}{{\Game^\ominus_n}}
\newcommand{\vare}{\varepsilon}

\begin{document}

\title{Chasing Silver}

\author{Andrzej Ros{\l}anowski}
\address{Department of Mathematics\\
 University of Nebraska at Omaha\\
 Omaha, NE 68182-0243, USA}

\email{roslanow@member.ams.org}

\urladdr{http://www.unomaha.edu/logic}

\author{Juris Stepr\={a}ns}
\address{Department of Mathematics\\
 York University\\
 4700 Keele Street\\
 Toronto, Ontario, Canada M3J 1P3}

\email{steprans@yorku.ca}
\urladdr{http://www.math.yorku.ca/$\sim$steprans}
\thanks{The research of the second author was partially supported by NSERC
  of Canada.}

\date{August 2005}

\begin{abstract}
Answering a question of the first author stated in \cite[0.2]{Ro0x} we show
that limits of CS iterations of $n$--Silver forcing notion have the
$n$--localization property. 
\end{abstract}

\maketitle

\section{Introduction}
The present paper is concerned with the $n$--localization property of the
$n$--Silver forcing notion and countable support (CS) iterations of such
forcings. The property of $n$--localization was introduced in Newelski and
Ros{\l}anowski \cite[p. 826]{NeRo93}.  
\begin{definition}
\label{nfor}
Let $n$ be an integer greater than 1. 
\begin{enumerate}
\item A tree $T$ is {\em an $n$--ary tree\/} provided that $(\forall s\in 
T)(|\suc_T(s)|\leq n)$. 
\item A forcing notion $\bbP$ has the {\em $n$--localization property\/} if
\[\forces_{\bbP}\mbox{`` }\big(\forall f\in\baire\big)\big(\exists T\in\bV
\big)\big(T\mbox{ is an $n$--ary tree and }f\in [T]\big)\mbox{ ''.}\]
\end{enumerate}
\end{definition}
Later the $n$--localization property, the $\sigma$--ideal generated by
$n$--ary trees and $n$--Sacks forcing notion $\bbD_n$ (see \ref{forcings})
have been found applied to problems on convexity numbers of closed subsets
of $\mbR^n$, see Geschke, Kojman, Kubi\'s and Schipperus \cite{GKKS02}, 
Geschke and Kojman \cite{GeKo02} and most recently Geschke \cite{Ge05}. 

We do not have any result of the form ``CS iteration of proper forcing
notions with the $n$--localization property has the $n$--localization''
yet. A somewhat uniform and general treatment of preserving the
$n$--localization has been recently presented in \cite{Ro0x}. However, the
treatment in that paper does not cover the $n$--Silver forcing notion
$\bbS_n$ (see \ref{forcings}), as a matter of fact it was not clear at some
moment if $\bbS_n$ has the property at all. It was stated in \cite[Theorem
2.3]{NeRo93} that the same proof as for $\bbD_n$ works also for CS
iterations and products of the $n$--Silver forcing notions $\bbS_n$ (see
Definition \ref{forcings}(3)). Maybe some old wisdom got lost, but it does
not look like that {\em the same arguments work for the $n$--Silver forcing
$\bbS_n$}. In the present paper we correct this gap and we provide a full
proof that CS iteartions of $\bbS_n$ (and other forcings listed in
\ref{forcings}) have the $n$--localization property, see \ref{seciter}. Our
main result \ref{main} on the $n$--Silver forcing seems to be very
$\bbS_n$--specific and it is not clear to which extend it may be
generalized. 

Let us explain what is a possible problem with the $n$--Silver forcing --
let us look at the ``classical'' Silver forcing $\bbS_2$. Given a Silver
condition $f$ such that $f\forces_{\bbS_2}\name{\tau}\in\baire$, standard
arguments allow it to be assumed that the complement of the domain of $f$
can be enumerated in the increasing order as $\{k_i:i<\omega\}$ and that for 
each $i\in\omega$ and $\rho:\{k_j: j<i\}\longrightarrow 2$ the condition
$f\cup \rho$ decides the value of $\name{\tau}\rest i$, say $f\cup\rho
\forces\name{\tau}\rest i=\sigma_\rho$. Now one could take the tree 
\[T^\oplus=\big\{\nu\in {}^{\omega{>}}\omega:\big(\exists i<\omega\big)
\big(\exists\rho\in {}^{\{k_j: j<i\}} 2\big)\big(\nu\trianglelefteq
\sigma_\rho\big)\big\}.\]
Easily $p\forces\name{\tau}\in [T^\oplus]$, but $T^\oplus$ does not have to
be a binary tree! (It could well be that $\sigma_\rho=\sigma^*$ for all
$\rho$ of length $100$ and then $\sigma_{\rho'}$ for $\rho'$ of length $101$
are pairwise distinct.) So we would like to make sure that $\sigma_\rho$ for 
$\rho$'s of the same length are distinct, but this does not have to be
possible. To show that $\bbS_2$ has the $2$--localization property we have
to be a little bit more careful. Let us give a combinatorial result which
easily implies that $\bbS_2$ has the $2$--localization property. Its proof
is in the heart of our proof of Theorem \ref{main}.  

Fix $\Psi:{}^{\omega{>}}2\longrightarrow\omega$. Define $\Psi^*:
{}^{\omega{>}} 2 \longrightarrow {}^{\omega{>}} \omega$ by induction: Let 
$\Psi^*(\langle\rangle) = \langle\rangle$ and define $\Psi^*(t\conc\langle
i\rangle)=\Psi^*(t)\conc\langle\Psi(t\conc\langle i\rangle)\rangle$. If
$\xi$ is a partial function from $\omega$ to $2$ and $\ell\leq\omega$
define  
\[W^\ell(\xi)=\big\{t \in {}^m 2:m< \min(\ell+1,\omega)\ \&\ \xi\rest m
\subseteq t\big\}\] 
and then define $T^\ell(\xi)=\{\Psi^*(t):t \in W^\ell(\xi)\}$, $T(\xi)=
T^\omega(\xi)$. 

\begin{theorem}
For any $\Psi:{}^{\omega{>}}2\longrightarrow\omega$ there is a partial
function $\xi:\omega\longrightarrow 2$ with co-infinite domain such that
$T(\xi)$ is a binary tree. 
\end{theorem}

\begin{proof}
To begin, two equivalence relations on ${}^{\omega{>}}2$ will be defined.  
First, define $s \equiv t$ if and only if $\Psi(t\conc \theta)=\Psi(s\conc
\theta)$ for all $\theta\in {}^{\omega{>}} 2$. Next, define $s\sim t$ if and
only if $\Psi^*(s)=\Psi^*(t)$.   

Now construct by induction on $m<\omega$ an increasing sequence
\[x_0< x_1 < \ldots < x_m < N_m\] 
and $\xi_m:N_m\setminus\{x_0,x_1\ldots x_m\}\longrightarrow 2$ such that 
$T^{N_m}(\xi_m)$ is a binary branching tree and, moreover, if $s$ and $t$
are maximal elements of $W^{N_m}(\xi_m)$ and $t\sim s$ then $t\equiv s$. The
induction starts with $N_0=0$. If the induction has been completed for $m$
then let $x_{m+1}=N_m$. Let $\Delta=\{d_0,d_1,\ldots,d_j\}$ be a set of
maximal elements of $T^{N_m}(\xi_m)$ such that precisely one member of each
$\sim$ equivalence class belongs to $\Delta$. Now, by induction on $i\leq j$
define $N^i$ and $\xi^i:N^i\setminus (N_m+1)\longrightarrow 2$ as
follows. Let $N^0=N_m+1$ and let $\xi^0=\emptyset$. Given $N^i$ and $\xi^i$,
if there is some $N> N^i$ and $\xi\supseteq \xi^i$ such that
$d_i\conc\langle 0\rangle\conc\xi\equiv d_i\conc\langle 1\rangle\conc\xi$
then let $N^{i+1}=N$ and let $\xi^{i+1}=\xi$. Otherwise it must be the case
that $d_i\conc\langle 0\rangle\conc \xi^i\not\equiv d_i\conc\langle
1\rangle\conc \xi^i$ and so it must be possible to find $N^{i+1}>N^i$ and
$\xi^{i+1}\supseteq\xi^i$ such that $d_i\conc\langle 0\rangle\conc\xi^{i+1} 
\not\sim d_i\conc\langle 1\rangle\xi^{i+1}$. Finally, let $N_{m+1}=N^j$ and
$\xi_{m+1}=\xi_m\cup\xi^j$. 

To see that this works, it must be shown that $T^{N_{m+1}}(\xi_{m+1})$ is  a
binary tree and that if $s$ and $t$ are maximal elements of
$W^{N_{m+1}}(\xi_{m+1})$ and $t\sim s$ then $t\equiv s$. To check the first
condition it suffices to take $t$ a maximal element of $T^{N_{m}}(\xi_{m})$
and check that the tree $T^{N_{m+1}}(\xi_{m+1})$ above $t$ is binary. Then
$t=\Psi^*(d_i)$ for some $i$ by the induction hypothesis. The tree
$T^{N_{m+1}}(\xi_{m+1})$ above $t$ consists is generated by all
$\Psi^*(d\conc\langle a\rangle\conc\xi_j)$ where $d\sim d_i$ and $a\in 2$. 
Note however that if $d\sim d_i$ then $d\equiv d_i$ and so
\[\Psi^*(d\conc\langle a\rangle\conc\xi_j)=\Psi^*(d_i\conc\langle
a\rangle\conc\xi_j).\]
Therefore $\Psi^*(d\conc\langle a\rangle\conc\xi_j)$ depends only on $a$ and 
not on $d$ and so $T^{N_{m+1}}(\xi_{m+1})$ is binary above $t$.

To check the second condition suppose that $s$ and $t$ are maximal elements
of $W^{N_{m+1}}(\xi_{m+1})$ and $t\sim s$. This implies that $t\rest N_m\sim
s\rest N_m$ and hence $t\rest N_m\equiv s\rest N_m$. Let $i$ be such that 
$t\rest N_m\sim s\rest N_m\sim d_i$. If $t(N_m)=s(N_m)=y$ then $t=t\rest
N_m\conc\langle y\rangle\conc\xi^j$ and $s= s\rest N_m\conc\langle
y\rangle\conc \xi^j$ and, since $t\rest N_m\equiv s\rest N_m$, it is
immediate that $t\equiv s$. So assume that $t(N_m)=0$ and $s(N_m)=1$. By the
same argument it follows that $t\equiv d_i\conc\langle 0\rangle\conc\xi^j$
and $s\equiv d_i\conc\langle 1\rangle\conc\xi^j$. Hence it suffices to show
that $d_i\conc\langle 0\rangle\conc\xi^j \equiv d_i\conc\langle
1\rangle\conc \xi^j$. Note that $d_i\conc\langle 0\rangle\conc\xi^j\sim
d_i\conc\langle 1\rangle\conc\xi^j$ since $t\sim d_i\conc\langle 0\rangle
\conc\xi^j$ and $s\sim d_i\conc\langle 1\rangle\conc\xi^j$. This means that
it must have been possible to find $\xi^i$ such that $d_i\conc\langle
0\rangle\conc\xi^i \equiv d_i\conc\langle 1\rangle\conc\xi^i$. It follows
that $d_i\conc\langle 0\rangle\conc\xi^j \equiv d_i\conc\langle 1\rangle
\conc\xi^j$. 
\end{proof}
\medskip

\noindent {\bf Notation:}\quad Our notation is rather standard and
compatible with that of classical textbooks (like Jech \cite{J}). In forcing
we keep the older convention that {\em a stronger condition is the larger
  one}.  

\begin{enumerate}
\item $n$ is our fixed integer, $n\geq 2$. 
\item For two sequences $\eta,\nu$ we write $\nu\vartriangleleft\eta$
whenever $\nu$ is a proper initial segment of $\eta$, and $\nu
\trianglelefteq\eta$ when either $\nu\vartriangleleft\eta$ or $\nu=\eta$. 
The length of a sequence $\eta$ is denoted by $\lh(\eta)$.
\item A {\em tree} is a family of finite sequences closed under initial
segments. For a tree $T$ and $\eta\in T$ we define {\em the successors of
$\eta$ in $T$} and {\em maximal points of $T$} by:
\[\begin{array}{rcl}
\suc_T(\eta)&=&\{\nu\in T: \eta\vartriangleleft\nu\ \&\ \neg(\exists\rho\in
T)(\eta\vartriangleleft\rho\vartriangleleft\nu)\},\\
\max(T)&=&\{\nu\in T:\mbox{ there is no }\rho\in T\mbox{ such that }
\nu\vartriangleleft\rho\}.
  \end{array}\]   
For a tree $T$ the family of all $\omega$--branches through $T$ is
denoted by $[T]$. 
\item For a forcing notion $\bbP$, all $\bbP$--names for objects in the
extension via $\bbP$ will be denoted with a tilde below (e.g.,
$\name{\tau}$, $\name{X}$). 
\end{enumerate}

\section{Definitions and the result}

\begin{definition}
\label{forcings}
\begin{enumerate}
\item {\em The $n$--Sacks forcing notion $\bbD_n$} consists of perfect trees  
  $p\subseteq {}^{\omega{>}}n$ such that 
\[(\forall \eta\in p)(\exists\nu\in p)(\eta\vtl\nu\ \&\ \suc_p(\eta)=n).\] 
The order of $\bbD_n$ is the reverse inclusion, i.e., $p\leq_{\bbD_n} q$ if
and only if $q\subseteq p$. (See \cite{NeRo93}.)

\item {\em The uniform $n$--Sacks forcing notion $\bbQ_n$} consists of
  perfect trees $p\subseteq {}^{\omega{>}}n$ such that 
\[(\exists X\in [\omega]^\omega)(\forall \eta\in p)(\lh(\eta)\in X\
  \Rightarrow\ \suc_p(\nu)=n).\]
The order of $\bbQ_n$ is the reverse inclusion, i.e., $p\leq_{\bbQ_n} q$ if
and only if $q\subseteq p$. (See \cite{Ro94}.)

\item Let us assume that $G=(V,E)$ is a hypergraph on a Polish space $V$
  such that 
\begin{itemize}
\item $E\subseteq [V]^{n+1}$ is open in the topology inherited from
  $V^{n+1}$, and  
\item $\big(\forall e\in E\big)\big(\forall v\in V\setminus e\big)\big(
  \exists w\in e\big)\big((e\setminus\{w\})\cup \{v\}\in E\big)$,   
\item for every non-empty open subset $U$ of $V$ and every countable family
  $\cF$ of subsets of $U$, either $\bigcup\cF\neq U$ or $[F]^{n+1}\cap E\neq
  \emptyset$ for some $F\in\cF$.  
\end{itemize}
{\em The Geschke forcing notion $\bbP_G$ for $G$} consists of all closed
sets $C\subseteq V$ such that the hypergraph $(C,E\cap [C]^{n+1})$ is
uncountably chromatic on every non-empty open subset of $C$. The order of
$\bbP_G$ is the inverse inclusion, i.e., $C\leq_{\bbP_G} D$ if and only  if
$D\subseteq C$. (See \cite{Ge05}.)
\end{enumerate}
\end{definition}

\begin{definition}
\label{silv}
\begin{enumerate}
\item {\em The $n$--Silver forcing notion $\bbS_n$} consists of partial
  functions $f$ such that $\Dom(f)\subseteq\omega$, $\Rng(f)\subseteq n$ 
  and $\omega\setminus\Dom(f)$ is infinite. The order of $\bbS_n$ is the
  inclusion, i.e., $f\leq_{\bbQ_n} g$ if and only  if $f\subseteq g$.
\item For an integer $i\in\omega$ and a condition $f\in\bbS_n$ we let
$\FP_i(f)$ to be the unique element of $\omega\setminus\Dom(f)$ such that  
$|\FP_i(f)\setminus \Dom(f)|=i$. (The $\FP$ stands for {\em Free Point}.) 
\item A binary relation $\leq^*_i$ on $\bbS_n$ is defined by 

$f\leq^*_i g$ if and only if ($f,g\in\bbS_n$ and) $f\leq_{\bbS_n} g$ and 
\[\big(\forall j\in\omega\big)\big(j<\lfloor i/4\rfloor\ \Rightarrow\
\FP_j(f)=\FP_j(g) \big).\]  
\item For $f\in\bbS_n$ and $\sigma:N\longrightarrow n$, $N<\omega$ we define
  $f*\sigma$ as the unique condition in $\bbS_n$ such that $\Dom(f*\sigma)=
  \Dom(f)\cup\{\FP_i(f):i<N\}$, $f\subseteq f*\sigma$ and $f*\sigma(\FP_i(f
  ))=\sigma(i)$ for $i<N$. 
\end{enumerate}
\end{definition}

\begin{definition}
\label{da2}
Let $\bbP$ be a forcing notion.
\begin{enumerate}
\item For a condition $p\in\bbP$ we define a game $\tGs(p,\bbP)$ of two
  players, {\em Generic\/} (``she'') and {\em Antigeneric\/} (``he''). A
  play of $\tGs(p,\bbP)$ lasts $\omega$ moves and during it the players
  construct a sequence  $\langle (s_i,\bar{\eta}^i,\bar{p}^i,\bar{q}^i):
  i<\omega\rangle$ as follows. At a stage $i<\omega$ of the play, first
  Generic chooses a finite $n$--ary tree $s_i$ such that 
\begin{enumerate}
\item[$(\alpha)$] $|\max(s_0)|\leq n$ and if $i=j+1$ then $s_j$ is a subtree
of $s_i$ such that  
\[\big(\forall \eta\in\max(s_i)\big)\big(\exists \ell<\lh(\eta)\big)\big(
\eta\rest \ell\in \max(s_j)\big),\] 
and 
\[\big(\forall\nu\in\max(s_j)\big)\big(0<\big|\big\{\eta\in\max(s_i):\nu
\vtl\eta\big\}\big|\leq n\big).\] 
\end{enumerate}
Next 
\begin{enumerate}
\item[$(\beta)$] Generic picks an enumeration $\bar{\eta}^i=\langle
  \eta^i_\ell:\ell<k_i\rangle$ of $\max(s_i)$ (so $k_i<\omega$) 
\end{enumerate}
and then the two players play a subgame of length $k_i$ choosing successive
terms of a sequence $\langle p^i_{ \eta^i_\ell},q^i_{\eta^i_\ell}:\ell<k_i
\rangle$. At a stage $\ell<k_i$ of the subgame, first Generic picks a
condition $p^i_{\eta^i_\ell}\in\bbP$ such that     
\begin{enumerate}
\item[$(\gamma)^i_\ell$] if $j<i$, $\nu\in \max(s_j)$ and $\nu\vtl
\eta^i_\ell$, then $q^j_\nu\leq p^i_{\eta^i_\ell}$ and $p\leq
p^i_{\eta^i_\ell}$,    
\end{enumerate}
and then Antigeneric answers with a condition $q^i_{\eta^i_\ell}$ stronger
than $p^i_{\eta^i_\ell}$. 

Finally, Generic wins the play $\langle (s_i,\bar{\eta}^i,\bar{p}^i,
\bar{q}^i):i<\omega\rangle$ if and only if 
\begin{enumerate}
\item[$(\circledast)$] there is a condition $q\geq p$ such that for every
$i<\omega$ the family $\{q^i_\eta:\eta\in\max(s_i)\}$ is predense above $q$.  
\end{enumerate}
\item We say that $\bbP$ has the {\em $\ominus_n$--property\/} whenever
Generic has a winning strategy in the game $\tGs(p,\bbP)$ for any
$p\in\bbP$. 
\item Let $K\in[\omega]^{\textstyle\omega}$, $p\in\bbP$. A strategy $\st$
  for Generic in $\tGs(p,\bbP)$ is {\em $K$--nice\/} whenever  
\begin{enumerate}
\item[$(\boxtimes_{\rm nice}^K)$] if so far Generic used $\st$ and $s_i$ and
  $\bar{\eta}^i=\langle \eta^i_\ell:\ell<k\rangle$ are given to her as
  innings at a stage $i<\omega$, then     
\begin{itemize}
\item $s_i\subseteq \bigcup\limits_{j\leq i+1} {}^j (n+1)$,
  $\max(s_i)\subseteq {}^{(i+1)}(n+1)$ and
\item if $\eta\in\max(s_i)$ and $i\notin K$, then $\eta(i)=n$, and 
\item if $\eta\in\max(s_i)$ and $i\in K$, then $\suc_{s_i}(\eta\rest i)=n$,
\item if $i\in K$ and $\langle p^i_{\eta^i_\ell},q^i_{\eta^i_\ell}:\ell<k
  \rangle$ is the result of the subgame of level $i$ in which Generic uses
  $\st$, then the conditions $p^i_{\eta^i_\ell}$ (for $\ell<k$) are pairwise
  incompatible. 
\end{itemize}
\end{enumerate}
\item We say that $\bbP$ has the {\em nice $\ominus_n$--property\/} if for
every $K\in [\omega]^{\textstyle\omega}$ and $p\in\bbP$, Generic has a
$K$--nice winning strategy in $\tGs(p,\bbP)$.
\end{enumerate}
\end{definition}

\begin{theorem}
[See {\cite[3.1+1.6+1.4]{Ro0x}}]
\label{prev}
The limits of CS iterations of the forcing notions defined in
\ref{forcings}, \ref{silv} have the nice $\ominus_n$--property.   
\end{theorem}

Now we may formulate our main result. 

\begin{theorem}
\label{main}
Assume that $\bbP$ has the nice $\ominus_n$--property and the
$n$--localization property. Let $\name{\bbS}_n$ be the $\bbP$--name for the
$n$--Silver forcing notion. Then the composition $\bbP*\name{\bbS}_n$ has
the $n$--localization property.  
\end{theorem}

The proof of Theorem \ref{main} is presented in the following section. Let
us note here that this theorem implies $n$--localization for CS iterations
of the forcing notions mentioned here. 

\begin{corollary}
\label{seciter}
Let $\bar{\bbQ}=\langle\bbP_\xi,\name{\bbQ}_\xi:\xi<\gamma\rangle$ be a
CS iteration such that, for every $\xi<\gamma$, $\name{\bbQ}_\xi$ is a
$\bbP_\xi$--name for one of the forcing notions defined in \ref{forcings},
\ref{silv}. Then $\bbP_{\gamma}=\lim(\bar{\bbQ})$ has the $n$--localization
property. 
\end{corollary}

\begin{proof}
By induction on $\gamma$. 

\noindent If $\gamma=\gamma_0+1$ and $\name{\bbQ}_{\gamma_0}$ is a
$\bbP_{\gamma_0}$--name for the $n$--Silver forcing notion, then
\ref{main} applies. (Note that $\bbP_{\gamma_0}$ has the nice
$\ominus_n$--property by \ref{prev} and it has the $n$--localization
property by the inductive hypothesis.) 

\noindent If $\gamma=\gamma_0+1$ and $\name{\bbQ}_{\gamma_0}$ is a
$\bbP_{\gamma_0}$--name for $\bbD_n$ or $\bbQ_n$ or $\bbP_G$, then
\cite[Theorem 3.4]{Ro0x} applies. (Note that $\bbP_{\gamma_0}$ has the
nice $\ominus_n$--property by \ref{prev} and it has the $n$--localization 
property by the inductive hypothesis.)   

\noindent If $\gamma$ is limit then \cite[3.5]{Ro0x} applies. 
\end{proof}

\begin{corollary}
No CS iteration of $\bbS_2$ adds an $\bbS_4$--generic real.
\end{corollary}

\section{Proof of Theorem {\protect\ref{main}}}

Let $\name{\tau}$ be a $\bbP*\name{\bbS}_n$--name for a member of
$\baire$. We may assume that $\forces_{\bbP*\name{\bbS}_n}\name{\tau}\notin
\bV$. If $G\subseteq\bbP$ is generic over $\bV$, then we will use the same
notation $\name{\tau}$ for $\bbS_n$--name in $\bV[G]$ for a member of
$\baire$ that is given by the original $\name{\tau}$ in the extension via
$\bbP*\name{\bbS}_n$. 

Let $(p,\name{f})\in \bbP*\name{\bbS}_n$ and let $\st$ be a winning
strategy of Generic in $\tGs(p,\bbP)$ which is nice for the set $K=\{4j+2:
j\in\omega\}$ (see \ref{da2}(3)).  

By induction on $i$ we are going to choose for each $i<\omega$
\[s_i,\bar{\eta}^i,\bar{p}^i,\bar{q}^i,\name{f}_i, \]
and for also $m_i,\bar{\sigma}^i$ for odd $i<\omega$ such that the following
conditions $(\boxtimes)_1$--$(\boxtimes)_7$ are satisfied.
\begin{enumerate}
\item[$(\boxtimes)_1$] $\langle s_i,\bar{\eta}^i,\bar{p}^i,\bar{q}^i:
  i<\omega\rangle$ is a play of $\tGs(p,\bbP)$ in which Generic uses $\st$. 
\item[$(\boxtimes)_2$] $\name{f}_i$ is a $\bbP$--name for a condition in
  $\bbS_n$, and we stipulate that $\name{f}_{-1}=\name{f}$. 
\item[$(\boxtimes)_3$] $q_\eta^i\forces_\bbP\name{f}_{i-1}\leq^*_i
  \name{f}_i$ for each $\eta\in\max(s_i)$. 
\end{enumerate}
For odd $i<\omega$:
\begin{enumerate}
\item[$(\boxtimes)_4$] $m_i<m_{i+2}<\omega$, $\bar{\sigma}^i=\langle
  \sigma^i_{\rho,\eta}:\eta\in \max(s_i)\ \&\ \rho\in {}^{\lfloor i/4
  \rfloor}n\rangle$, $\sigma^i_{\rho,\eta}:m_i\longrightarrow\omega$. 
\item[$(\boxtimes)_5$] $(q_\eta^i,\name{f}_i*\rho)\forces_{\bbP*
  \name{\bbS}_n}$`` $\name{\tau}\rest m_i=\sigma^i_{\rho,\eta}$ '' for
  $\rho\in {}^{\lfloor i/4\rfloor}n$ and $\eta\in\max(s_i)$.  
\item[$(\boxtimes)_6$] If $\eta\in\max(s_i)$ and $\rho,\rho':\lfloor i/4
  \rfloor\longrightarrow n$ are distinct but $\sigma^i_{\rho,\eta}=
  \sigma^i_{\rho',\eta}$, then for every $q\geq q_\eta^i$ and a $\bbP$--name
  $\name{g}$ for an $n$--Silver condition and $m,\sigma,\sigma'$ such that  
\[q\forces_{\bbP}\name{f}_i\leq^*_i\name{g},\quad (q,\name{g}*\rho)
  \forces_{\bbP*\name{\bbS}_n}\name{\tau}\rest m=\sigma,\quad (q,\name{g}*
  \rho')\forces_{\bbP*\name{\bbS}_n}\name{\tau}\rest m=\sigma'\] 
we have $\sigma=\sigma'$.
\item[$(\boxtimes)_7$] If $\eta,\eta'\in\max(s_i)$ are distinct, $\rho,
  \rho':\lfloor i/4\rfloor\longrightarrow n$, then $\sigma^i_{\rho,\eta}\neq
  \sigma^i_{\rho',\eta'}$.
\end{enumerate}
\medskip

So suppose that $i<\omega$ is even and we have already defined $s_{i-1},
\bar{q}^{i-1},m_{i-1}$ and $\name{f}_{i-1}$ (we stipulate $s_{-1}=
\{\langle\rangle\}$, $q^{-1}_{\langle\rangle}=p$, $\name{f}_{-1}=\name{f}$
and $m_{-1}=0$). Let $j=\lfloor i/4\rfloor$ (so either $i=4j$ or $i=4j+2$). 
\medskip

The strategy $\st$ and demand $(\boxtimes)_1$ determine $s_i$ and
$\bar{\eta}^i=\langle\eta^i_k:k<k_i\rangle$. To define $\bar{p}^i,\bar{q}^i$
and $\name{f}_i$ we consider the following run of the subgame of level $i$
of $\tGs(p,\bbP)$. Assume we are at stage $k<k_i$ of the subgame. Now,
$p^i_{\eta^i_k}$ is given by the strategy $\st$ (and $(\boxtimes)_1$, of 
course). Suppose for a moment that $G\subseteq\bbP$ is generic over $\bV$,
$p^i_{\eta^i_k}\in G$. Working in $\bV[G]$ we may choose $\bar{\ell},
\bar{L},g^*,\bar{\sigma}^*, M$ such that     
\begin{enumerate}
\item[$(\boxtimes)_8^\alpha$] $M=n^j$, $\bar{\ell}=\langle\ell_m:m\leq
  M\rangle$ and $j=\ell_0<\ldots<\ell_M$, $\bar{L}=\langle L_m:m\leq
  M\rangle$ and $m_{i-1}<L_0<\ldots<L_M$,   
\item[$(\boxtimes)_8^\beta$]  $g^*\in\bbS_n$, $\name{f}_{i-1}[G]\leq^*_i
  g^*$ and $\bar{\sigma}^*=\langle\sigma^*_\rho:\rho\in {}^{\ell_M}
  n\rangle$, $\sigma^*_\rho\in {}^{L_M}\omega$ (for $\rho\in {}^{\ell_M}n$),  
\item[$(\boxtimes)_8^\gamma$] $g^**(\rho\rest\ell_m)\forces_{\bbS_n}$``
  $\name{\tau}\rest L_m= \sigma^*_\rho\rest L_m$ '' for each $m\leq M$ and
  $\rho\in {}^{\ell_M}n$,    
\item[$(\boxtimes)_8^\delta$] if $\rho_0,\rho_1\in {}^{\ell_M} n$,
  $\rho_0\rest j\neq \rho_1\rest j$ but $\sigma^*_{\rho_0}\rest
  L_0=\sigma^*_{\rho_1}\rest L_0$, then there is no condition $g\in\bbS_n$
  such that $g^*\leq^*_i g$ and for some $L<\omega$ and distinct $\sigma_0,
  \sigma_1\in {}^L\omega$ we have that $g*\rho_0\forces\name{\tau}\rest L
  =\sigma_0$, $g*\rho_1 \forces\name{\tau}\rest L=\sigma_1$,    
\item[$(\boxtimes)_8^\vare$]  for each $m<M$ and $\rho_0\in {}^{\ell_m}n$
  the set $\{\sigma^*_\rho\rest [L_m,L_{m+1}):\rho_0\vtl\rho\in {}^{\ell_M}
  n\}$ has at least $n^j\cdot k_i+777$ elements. 
 \end{enumerate}
It should be clear how the construction is done. (First we take care of
clause $(\boxtimes)_8^\delta$ by going successively through all pairs of 
elements of ${}^jn$ and trying to force distinct values for initial segments 
of $\name{\tau}$, if only this is possible. Then we ensure
$(\boxtimes)_8^\vare$ basically by deciding longer and longer initial
segments of $\name{\tau}$ on fronts/levels of a fusion sequence of
conditions in $\bbS_n$ and using the assumption that $\name{\tau}$ is forced
to be ``new''.) Now, going back to $\bV$, we may choose a condition
$q^i_{\eta^i_k}\in\bbP$ stronger than $p^i_{\eta^i_k}$ and a $\bbP$--name
$\name{g}^{*,k}$ for a condition in $\bbS_n$ and objects $\bar{\ell}^k,
\bar{L}^k,\bar{\sigma}^{*,k}$ such that   
\[q^i_{\eta^i_k}\forces_{\bbP}\mbox{`` }\bar{\ell}^k,\bar{L}^k,
\name{g}^{*,k},\bar{\sigma}^{*,k},n^j\mbox{ satisfy clauses
  $(\boxtimes)_8^\alpha$--$(\boxtimes)_8^\vare$ as $\bar{\ell},\bar{L},g^*, 
  \bar{\sigma}^*,M$ there ''.}\] 
The condition $q^i_{\eta^i_k}$ is treated as an inning of Antigeneric at
stage $k$ of the subgame of $\tGs(p,\bbP)$ and the process continues.  

After the subgame of level $i$ is completed, we have defined $\bar{p}^i$ and
$\bar{q}^i$. We also choose $\name{f}_i$ to be a $\bbP$--name for an element
of $\name{\bbS}_n$ such that $\forces_\bbP$`` $\name{f}_{i-1}\leq^*_i
\name{f}_i$ '' and $q^i_{\eta^i_k}\forces_\bbP$`` $\name{f}_i=\name{g}^{*,
k}$ '' for all $k<k_i$ (remember that $\st$ is nice, so the conditions
$q^i_{\eta^i_k}$ are pairwise incompatible). This completes the description
of what happens at the stage $i$ of the construction (one easily verifies
that $(\boxtimes)_1$--$(\boxtimes)_3$ are satisfied) and we proceed to the
next, $i+1$, stage. Note that $\lfloor (i+1)/4\rfloor=j$. 
\medskip

We let $m_{i+1}=\max(L^k_M: k<k_i)+5$ and let $\ell=\max(\ell^k_M:
k<k_i)+5$. Similarly as at stage $i$, $s_{i+1}$ and $\bar{\eta}^{i+1}=
\langle\eta^{i+1}_k:k<k_{i+1}\rangle$ are determined by the strategy $\st$
and $(\boxtimes)_1$; note that $\max(s_{i+1})=\{\nu\conc\langle n\rangle:
\nu\in\max(s_i)\}$ so $k_{i+1}=k_i$. To define $\bar{p}^{i+1},\bar{q}^{i+1}$
and $\name{f}_{i+1}$ we consider the following round of the subgame of level
$i+1$ of $\tGs(p,\bbP)$. At a stage $k<k_{i+1}$ of the subgame, letting
$\eta=\eta^{i+1}_k$, the condition $p^{i+1}_\eta$ is given by the strategy
$\st$. Suppose for a moment that $G\subseteq\bbP$ is generic over $\bV$,
$p^{i+1}_\eta\in G$. In $\bV[G]$ we may choose a condition $h^*\in \bbS_n$
such that 
\begin{enumerate}
\item[$(\boxtimes)_9$] $\name{f}_i[G]\leq^*_\ell h^*$ and for every
$\rho\in {}^\ell n$ the condition $h^**\rho$ decides the value of
$\name{\tau}\rest m_{i+1}$, say $h^**\rho\forces_{\bbS_n}$`` $\name{\tau}  
\rest m_{i+1}=\sigma_\rho$ ''. 
\end{enumerate}
Then going back to $\bV$ we choose a $\bbP$--name $\name{h}^{*,\eta}$ for a
condition in $\bbS_n$, a sequence $\bar{\sigma}^\eta=\langle
\sigma^\eta_\rho: \rho\in {}^\ell n\rangle$ and a condition
$q^{i+1}_\eta\geq p^{i+1}_\eta$ such that    
\[q^{i+1}_\eta\forces_\bbP\mbox{`` }\name{h}^{*,\eta},\bar{\sigma}^\eta
\mbox{ are as in $(\boxtimes)_9$ ''.}\] 
The condition $q^{i+1}_\eta$ is treated as an inning of Antigeneric at stage
$k$ of the subgame of $\tGs(p,\bbP)$ and the process continues.   

After the subgame of level $i+1$ is completed, we have defined
$\bar{p}^{i+1}$ and $\bar{q}^{i+1}$. Since for every $\eta\in\max(s_{i+1})$
we have that $p^{i+1}_\eta\geq q^i_{\eta\rest (i+1)}$, we may use
$(\boxtimes)^\vare_8$ and choose $\rho(\eta):[j,\ell)\longrightarrow n$ (for
$\eta\in\max(s_{i+1})$) such that   
\begin{enumerate}
\item[$(\boxtimes)_{10}$] if $\eta,\eta'\in\max(s_{i+1})$ are distinct and
$\theta,\theta'\in {}^j n$, and $\rho=\theta\conc\rho(\eta)$, $\rho'=\theta'
\conc\rho(\eta')$, then $\sigma^\eta_\rho\neq\sigma^{\eta'}_{\rho'}$.   
\end{enumerate}
Let $\name{f}_i$ be a $\bbP$--name for a condition in $\bbS_n$ such that
$\forces_{\bbP} \name{f}_i\leq^*_{i+1}\name{f}_{i+1}$ and 
\[q^{i+1}_\eta\forces_{\bbP}\mbox{`` }\name{h}^{*,\eta}\leq^*_i\name{f}_i\
\&\ \big(\forall\theta\in {}^jn\big)\big(\name{f}_i*\theta=\name{h}^{*,\eta}
*(\theta\conc \rho(\eta))\big)\mbox{ ''.}\] 
Also, for $\eta\in\max(s_{i+1})$ and $\rho\in {}^jn$, we let $\sigma^{i+
1}_{\rho,\eta}=\sigma^\eta_{\rho\conc^\rho(\eta)}$. This completes the
description of what happens at the stage $i+1$ of the construction (one
easily checks that $(\boxtimes)_1$--$(\boxtimes)_7$ are satisfied). Thus we
have finished the description of the inductive step of the construction of
$s_i,\bar{\eta}^i,\bar{p}^i,\bar{q}^i,\name{f}_i$ (for $i<\omega$). 
\bigskip

After the construction is carried out we may pick a condition $q\in\bbP$
stronger than $p$ and such that for each $i<\omega$ the family
$\{q^i_\eta:\eta\in\max(s_i)\}$ is predense above $q$ (possible by
$(\boxtimes)_1$). 

Suppose that $G\subseteq\bbP$ is generic over $\bV$, $q\in G$. Then there is
$\eta\in {}^\omega (n+1)$ such that $\eta\rest (i+1)\in \max(s_i)$ and
$q^i_{\eta\rest (i+1)}\in G$ for each $i<\omega$. Therefore we may use
$(\boxtimes)_3$ to conclude that there is a condition $g\in \bbS_n$ stronger
than all $\name{f}_i[G]$. Going back to $\bV$, we may choose a $\bbP$--name
$\name{g}$ for a condition in $\bbS_n$ such that $q\forces_\bbP(\forall
i<\omega)(\name{f}_i\leq \name{g})$. 

Note that for each $i<\omega$ the family $\big\{(q^i_\eta,\name{f}_i*\rho):
\eta\in\max(s_i)\ \&\ \rho\in {}^{\lfloor i/4\rfloor} n\big\}$
is predense in $\bbP*\name{\bbS}_n$ above $(q,\name{g})$, and hence (by
$(\boxtimes)_5$) 
\[(q,\name{g})\forces_{\bbP*\name{\bbS}_n}\mbox{`` }\name{\tau}\rest m_i\in
\{\sigma^i_{\rho,\eta}:\eta\in\max(s_i)\ \&\ \rho\in {}^{\lfloor i/4\rfloor}
n\}\mbox{ for every odd $i<\omega$ ''.}\]
Also,  
\begin{enumerate}
\item[$(\boxtimes)_{11}$] if $i\geq 3$ is odd, $\eta\in\max(s_i)$, $\rho
  \in {}^{\lfloor i/4\rfloor} n$ and $\eta'=\eta\rest (i-1)$ and $\rho'=
  \rho\rest\lfloor (i-2)/4\rfloor$, then $\eta'\in\max(s_{i-2})$ and
  $\sigma^{i-2}_{\rho',\eta'}=\sigma^i_{\rho,\eta}\rest m_{i-2}$. 
\end{enumerate}
[Why? Since $\st$ is a nice strategy, $\eta\rest i\in\max(s_{i-1})$ and
  $\eta'\in\max(s_{i-2})$. It follows from $(\boxtimes)_1$ that
  $q^{i-2}_{\eta'}\leq q^{i-1}_{\eta\rest i}\leq q^i_\eta$ and by
  $(\boxtimes)_3$ we have $q^i_\eta\forces_\bbP \name{f}_{i-2}
  \leq^*_{i-1}\name{f}_i$. Therefore $q^i_\eta\forces_\bbP
  \name{f}_{i-2}*\rho'\leq\name{f}_i*\rho$ and $(q^{i-2}_{\eta'}, 
  \name{f}_{i-2}*\rho')\leq (q^i_\eta,\name{f}_i*\rho)$, so using
  $(\boxtimes)_5$ we may conclude that $\sigma^{i-2}_{\rho',\eta'}=
  \sigma^i_{\rho,\eta}\rest m_{i-2}$.]\\   
Let 
\[T=\big\{\nu\in {}^{\omega{>}}\omega: \big(\exists i<\omega\mbox{ odd\/}
\big)\big(\exists\eta\in\max(s_i)\big)\big(\exists\rho\in {}^{\lfloor i/4
\rfloor}n \big)\big(\nu\trianglelefteq \sigma^i_{\rho,\eta}\big)\big\}.\]  
Then $T$ is a perfect tree and $(q,\name{g})\forces_{\bbP*\name{\bbS}_n}
\name{\tau}\in [T]$. So the theorem will readily follow once we show that
$T$ is $n$---ary. To this end we are going to argue that 
\begin{enumerate}
\item[$(\boxtimes)_{12}$] if $i\geq 3$ is odd, $\eta\in\max(s_i)$, $\rho\in
  {}^{\lfloor i/4\rfloor} n$, then 
\[\big|\big\{\sigma^i_{\pi,\nu}: \nu\in \max(s_i)\ \ \&\ \ \pi\in
{}^{\lfloor i/4\rfloor} n\ \ \&\ \ \sigma^i_{\rho,\eta}\rest m_{i-2}=  
\sigma^i_{\pi,\nu}\rest m_{i-2}\big\}\big|\leq n.\]   
\end{enumerate}
\medskip

\noindent {\sc Case A:}\qquad $i=4j+1$ for some $j<\omega$.\\ 
Suppose that $\eta,\nu\in\max(s_i)$, $\rho,\pi\in {}^{\lfloor i/4\rfloor}
n$ are such that $\sigma^i_{\rho,\eta}\neq\sigma^i_{\pi,\nu}$ but
$\sigma^i_{\rho,\eta}\rest m_{i-2}=\sigma^i_{\pi,\nu}\rest m_{i-2}$. The
latter and $(\boxtimes)_7$ imply that $\eta\rest (i-1)=\nu\rest (i-1)$, and
since $i-1,i\notin K$ we get that $\eta(i-1)=\nu(i-1)=n=\eta(i)=\nu(i)$
(remember: $\st$ is nice for $K$), so $\eta=\nu$. If $\rho\rest (j-1)\neq
\pi\rest (j-1)$, then let $\rho'=\rho\rest(j-1)\conc\langle\pi(j-1)\rangle$,
otherwise $\rho'=\pi$. 

Suppose $\rho'\neq\pi$. Let $\name{g}$ be (a $\bbP$--name for) $\name{f}_i
\cup\big\{\big(\FP_{j-1}(\name{f}_i),\pi(j-1)\big)\big\}$ and
$q=q^i_\eta$. Then $q\geq q^{i-2}_{\eta\rest (i-1)}$, $q\forces
\name{f}_{i-2}\leq^*_{i-2}\name{g}$, and 
\[q\forces\mbox{`` }\name{g}*\big(\rho'\rest(j-1)\big)=\name{f}_i*\rho'
\mbox{ and }\name{g}*\big(\pi\rest(j-1)\big)=\name{f}_i*\pi\mbox{ ''.}\]
Hence 
\[(q,\name{g}*\big(\rho'\rest(j-1)\big)\forces\mbox{`` }\name{\tau}\rest
m_i=\sigma^i_{\rho',\eta}\mbox{ ''\quad and\quad }(q,\name{g}*\big(\pi
\rest(j-1)\big)\forces\mbox{`` }\name{\tau}\rest m_i=\sigma^i_{\pi,\eta}
\mbox{ ''.}\]
Now we use our assumption that $\sigma^i_{\rho',\eta}\rest m_{i-2}=
\sigma^i_{\pi,\eta}\rest m_{i-2}$ (and $(\boxtimes)_{11}$) and
$(\boxtimes)_6$ to conclude that $\sigma^i_{\rho',\eta}=\sigma^i_{\pi,
  \eta}$. Trivially the same conclusion holds if $\rho'=\pi$, so we have
justified that 
\[\begin{array}{r}
\big\{\sigma^i_{\pi,\nu}: \nu\in \max(s_i)\ \ \&\ \ \pi\in {}^{\lfloor
i/4\rfloor} n\ \ \&\ \ \sigma^i_{\rho,\eta}\rest m_{i-2}=\sigma^i_{\pi,
  \nu}\rest m_{i-2}\big\}\subseteq\qquad\\ 
\big\{\sigma^i_{\pi,\eta}: \pi\in {}^j n\ \ \&\ \ \rho\rest (j-1)=\pi\rest
(j-1)\big\}
  \end{array}\]
and the latter set is of size at most $n$. 
\medskip

\noindent {\sc Case B:}\qquad $i=4j+3$ for some $j<\omega$.\\ 
Again, let us assume that $\eta,\nu\in\max(s_i)$, $\rho,\pi\in {}^{\lfloor
  i/4\rfloor} n$ are such that $\sigma^i_{\rho,\eta}\neq\sigma^i_{\pi,\nu}$
but $\sigma^i_{\rho,\eta}\rest m_{i-2}=\sigma^i_{\pi,\nu}\rest m_{i-2}$. Then,
like in the previous case, $(\boxtimes)_7$ implies $\eta\rest (i-1)=\nu\rest
(i-1)$. Also $\lfloor i/4\rfloor=j=\lfloor (i-2)/4\rfloor$, so $\rho\rest
\lfloor (i-2)/4\rfloor=\rho$, $\pi\rest \lfloor (i-2)/4\rfloor=\pi$. Now, if
$\rho=\pi$, then trivially $\sigma^i_{\pi,\nu}=\sigma^i_{\rho,\nu}$. If
$\rho\neq\pi$, then we use $(\boxtimes)_6$ (with $i-2,\rho,\pi,q^i_\eta,
\name{f}_i$ here in place of $i,\rho,\rho',q,\name{g}$ there, respectively)
to argue that $\sigma^i_{\pi,\nu}=\sigma^i_{\rho,\nu}$. Consequently 
\[\begin{array}{r}
\big\{\sigma^i_{\pi,\nu}: \nu\in\max(s_i)\ \&\ \pi\in {}^{\lfloor i/4
\rfloor} n\ \&\ \sigma^i_{\rho,\eta}\rest m_{i-2}=\sigma^i_{\pi,\nu}\rest
m_{i-2}\big\}\subseteq\qquad\\ 
\big\{\sigma^i_{\rho,\nu}:\nu\in\max(s_i)\ \&\ \eta\rest (i-2)=\nu\rest
(i-2)\big\} 
  \end{array}\]
and the latter set is of size at most $n$. 

Now in both cases we easily get the assertion of $(\boxtimes)_{12}$,
completing the proof of the theorem.


\begin{thebibliography}{1}
\bibitem{Ge05}
Stefan Geschke.
\newblock {More on convexity numbers of closed sets in ${\mathbb R}^n$}.
\newblock {\em Proceedings of the American Mathematical Society},
  133:1307--1315, 2005.

\bibitem{GeKo02}
Stefan Geschke and Menachem Kojman.
\newblock {Convexity numbers of closed sets in ${\mathbb R}^n$}.
\newblock {\em Proceedings of the American Mathematical Society},
  130:2871--2881, 2002.

\bibitem{GKKS02}
Stefan Geschke, Menachem Kojman, Wies{\l}aw Kubi\'s, and Rene Schipperus.
\newblock {Convex decompositions in the plane and continuous pair colorings of
  the irrationals}.
\newblock {\em Israel Journal of Mathematics}, 131:285--317, 2002.

\bibitem{J}
Thomas Jech.
\newblock {\em {Set theory}}.
\newblock Academic Press, New York, 1978.

\bibitem{NeRo93}
Ludomir Newelski and Andrzej Ros{\l}anowski.
\newblock {The ideal determined by the unsymmetric game}.
\newblock {\em Proceedings of the American Mathematical Society}, 117:823--831,
  1993.

\bibitem{Ro0x}
Andrzej Ros{\l}anowski.
\newblock {$n$--localization property}.
\newblock {\em Preprint}.
\newblock math.LO/0507519.

\bibitem{Ro94}
Andrzej Ros{\l}anowski.
\newblock Mycielski ideals generated by uncountable systems.
\newblock {\em Colloquium Mathematicum}, LXVI:187--200, 1994.
\end{thebibliography}


\end{document}